\documentclass[11pt]{amsart}
\usepackage{amsmath}
\usepackage{amssymb}

\newtheorem{theorem}{Theorem}[section]

\newtheorem{rgchtheorem}[theorem]{RGCH Theorem (partial version)}
\newtheorem{claim}[theorem]{Claim}

\newtheorem{lemma}[theorem]{Lemma}

\newtheorem{corollary}[theorem]{Corollary}

\newtheorem{fact}[theorem]{Fact}

\theoremstyle{definition}
\newtheorem{definition}[theorem]{Definition}

\theoremstyle{remark}
\newtheorem{remark}[theorem]{Remark}

\newcount\skewfactor
\def\mathunderaccent#1#2 {\let\theaccent#1\skewfactor#2
\mathpalette\putaccentunder}
\def\putaccentunder#1#2{\oalign{$#1#2$\crcr\hidewidth
\vbox
to.2ex{\hbox{$#1\skew\skewfactor\theaccent{}$}\vss}\hidewidth}}

\def\smallbox#1{\leavevmode\thinspace\hbox{\vrule\vtop{\vbox
   {\hrule\kern1pt\hbox{\vphantom{\tt/}\thinspace{\tt#1}\thinspace}}
   \kern1pt\hrule}\vrule}\thinspace}

\newcommand\includedin\subseteq
\newcommand\intersect\cap
\newcommand\union\cup

\newcommand{\cf}{{\rm cf}}
\newcommand{\Dom}{{\rm Dom}}

\renewcommand{\int}{{\rm int}}
\newcommand{\Min}{{\rm Min}}
\newcommand{\Sup}{{\rm Sup}}
\newcommand{\supRang}{{\rm supRang}}
\newcommand{\Rang}{{\rm Rang}}

\newcommand\GG{{\mathbb G}}

\newcommand\niceF{{\mathcal F}}
\newcommand\niceG{{\mathcal G}}
\newcommand\niceW{{\mathcal W}}
\newcommand\niceB{{\mathcal B}}
\newcommand\niceT{{\mathcal T}}
\newcommand\niceH{{\mathcal H}}
\newcommand\niceP{{\mathcal P}}

\newcommand\xfrak{{\mathfrak x}}
\newcommand\Mfrak{{\mathfrak M}}

\newcommand\ibold{{\mathbf i}}
\newcommand\jbold{{\mathbf j}}

\newcommand\gbold{{\mathbf g}}
\newcommand\hbold{{\mathbf h}}

\newcommand{\eps}{{\varepsilon}}

\newcommand{\st}{{such that }}

\newcommand{\rest}{{\restriction}}

\newcommand{\mat}{\mathcal}

\begin{document}
\title{ Existence of EF-equivalent Non Isomorphic Models} 
\author{Chanoch Havlin}
\author{Saharon Shelah}
\address{Institute of Mathematics\\
 The Hebrew University of Jerusalem\\
 Jerusalem 91904, Israel\\
 and\\  Department of Mathematics\\
 Rutgers University\\
 New Brunswick, NJ 08854, USA}
\email{shelah@math.huji.ac.il}
\urladdr{http://www.math.rutgers.edu/\char`\~shelah}
\thanks{We would like to thank the Israel 
Science Foundation for partial support of this research (Grant No. 242/03). 
Publication 866}

\begin{abstract}
We prove the existence of pairs of models of the same cardinality $\lambda$ which are very equivalent 
according to EF games, but not isomorphic. We continue the paper \cite{SH-836}, but we don't rely on it.
\end{abstract}

\maketitle

\tableofcontents
\newpage
%
%
\section{Introduction}
\noindent 
There had been much study of equivalence relations between models. When we study such an equivalence relation,
one of the basic questions is - is this relation actually trivial - equivalent models are isomorphic? For example, countable models
which are elementary equivalent in $L_{\omega_1,\omega}$ are isomorphic.\\ (Scott showed this in \cite{Scott} for countable vocabulary, and Chang generalized it in \cite{Chang} for any vocabulary). For $\lambda = \cf(\lambda) > \aleph_0$, Morely gave (without publishing) a counter example - a pair of $L_{\infty,\lambda}$ equivalent models of size $\lambda$ which are not isomorphic. Shelah(\cite{SH-Classification}) gave such an example for almost every singular $\lambda$ .\\
Those questions also relate to classification theory : The existence of "strongly" equivalent models which are not isomorphic is a non-structure property for a class of models. On the other side, if "not too strong" equivalence relation is actually the isomorphism relation, this is a structure property.(See \cite{SH-Classification} and
\cite{Hyttinen}).
\par \noindent\\
One of the equivalence relations studied in this context, is equivalence \\under EF( Ehernfeucht-Fraisse) games.
A detailed discussion of EF games and their history, can be found in \cite{Hodges}. The general structure of an EF game 
on a pair of models is as follows:\\
There are two players - isomorphism player, who we call ISO and anti-isomorphism player, who we call AIS. During the game, AIS chooses members of the models, and ISO defines "interactively" a partial isomorphism between the models
 - in every move he has to extend that partial isomorphism, such that the elements chosen by AIS will be contained in the domain or in the range. The isomorphism player looses the game if at some point, he cannot find a legal move. If he doesn't loose, he wins. We limit the length of the game and the number of elements that AIS may choose at each move.
(Because, if AIS can list all the members of one of the models, then the game is not interesting). In \cite{SH-836}, the games were with fixed length. In this paper, we deal with EF games approximated by trees - the length of the game is limited by adding the demand that in each move, AIS has to choose a node in some fixed tree $\niceT$ (with certain properties), such that the sequence of nodes formed by his choices, is strictly increasing in the order $<^\niceT$. If AIS cannot choose such node - he looses.\\ 
We say that two models are equivalent with respect to some EF game $\Game$ , if ISO has a winning strategy in $\Game$
played on those models.\\

\par\noindent 
In \cite{SH-836} it was proved that if $\lambda = \cf(\lambda) = \lambda^{\aleph_0}$, then there are non isomorphic models of size $\lambda$ which are $EF_{\alpha,\lambda}$ equivalent for every $\alpha < \lambda$. Where $EF_{\alpha,\lambda}$ equivalence  means that they are equivalent under every EF game with $\alpha$ stages, such that AIS has to choose $ < \lambda$ members of the models at each stage. There was also a result for $\lambda$ singular,  with a necessary change of the equivalence relation.\\
\par\noindent 
Here we generalize the results in 2 ways:
first, we move to EF games \\ approximated by trees instead of fixed-length games( See Hyttinen and\\ Tuuri in \cite{Hyttinen} who investigated such games in the context of classification\\ theory). Second, we give results also for $\lambda >\beth_\omega$ without the assumption $\lambda = \lambda^{\aleph_0}$, where we use PCF theory to have some "approximation" instead of $\lambda = \lambda^{\aleph_0}$.\\
\par \noindent
In section 1 we prove that for regular $\lambda = \lambda^{\aleph_{0}}$ for some class of reasonably large trees ( See detailed discussion justifying the choice, in the beginning of section 1 ) for every tree from that class there are non isomorphic models of size $\lambda$ which are equivalent under EF games approximated by that tree, such that in each move AIS is allowed to choose $ < \lambda$ members of the models ( See definition \ref{Sec1-Def-Game-With-Tree} ).\\
\par\noindent
In section 2 we do the parallel for singular $\lambda $. But for singular $\lambda$, if we allow AIS to choose $< \lambda$ elements in each move, and the tree has a branch of length $\cf(\lambda)$, then the game is not interesting, because AIS can choose all the members of the models during the game. So we have to be more careful - we allow AIS to choose only one element in each move. This is still a generalization of the result for such $\lambda$ in \cite{SH-836} - see the discussion at the beginning of section 2.\\
\par\noindent
In section 3 we prove that for regular $\lambda > \beth_\omega$, for every tree of size $\lambda$ without a branch of length $\lambda$, there are non isomorphic models of size $\lambda$ which are equivalent under the EF game approximated by that tree, such that in each move AIS is allowed to choose $ < \lambda$ members of the models.\\
\par\noindent
In section 4 we prove a similar result for $\lambda > \cf(\lambda) > \beth_\omega$. As we explained above, because of the singularity of $\lambda$, we have to restrict the number of elements that AIS is allowed to choose at each move - in stage $\alpha$, AIS has to choose $ < 1 + \alpha$ members of the models.

\newpage
%
%
\section{Games with trees for regular $\lambda = \lambda^{\aleph_{0}}$}
\noindent
In \cite{Hyttinen} there is a construction of non-isomorphic models of size $\lambda$ which are equivalent under EF
games approximated by trees of size $\lambda$ with no $\lambda$ branch, when $\lambda = \lambda^{<\lambda}$.
In \cite{SH-836} there is such a construction under a weaker assumption on $\lambda$ - $\lambda = \cf(\lambda) =
\lambda^{\aleph_0}$, but there the result is for games of any fixed \mbox{length $< \lambda$,} not for games
which are approximated by trees. We want to generalize this result to games approximated by trees.\\
Now, which trees should we consider ? If we limit ourselves only to trees of size $\lambda$, It seems that
the set of trees will be "small". Why? - assume for example that $\lambda= \cf(\lambda) = \lambda^{\aleph_0}
< \lambda^{\aleph_1}$. A tree of size $\lambda$ must drop at least one of the following conditions :
\begin{enumerate}
\item above every node there is an antichain of size $\lambda$
\item every chain of size $\leq \aleph_1$ has an upper bound
\end{enumerate}
If $\lambda \gg \aleph_1$, this kind of trees seem to be too
degenerate. We could have demanded that the size of the tree will
be $\leq 2^{< \lambda}$. But it is possible that $2^{< \lambda} =
2^\lambda$ and it is
reasonable to assume that the result will not be true in this case.\\
We take the middle road : we don't limit explicitly the size of
the tree, but we demand that the tree will be "definable" enough -
the cause of not having a branch of length $\lambda$, is that the
nodes of the tree are actually partial functions from $\lambda$ to
$\lambda$ which satisfy a certain local condition. By "local" we
mean that a function $f$ satisfies the condition iff any
restriction of the $f$ to a countable set satisfies it. The tree
order is inclusion, and there is no function from $\lambda$ to
$\lambda$ which satisfies the condition. By
\ref{Sec1-Remark-Strength-Of-Game} this result is indeed
generalization of "for every tree of size $\lambda$ and no
$\lambda$ branch".

\begin{definition}
\label{Sec1-Def-Game-With-Tree} For a tree $\niceT$, a cardinal
$\mu$, and models with common vocabulary $M_1,M_2$,\hspace{1mm}
the game $\Game_{\niceT,\mu}(M_1,M_2)$ between the players ISO and
AIS is defined as follows:\\
After stage $\alpha$ in the game we have the sequence $\langle
f_\beta :\,\beta \leq \alpha\,\rangle$, which is an increasing
continuous sequence of partial isomorphisms from $M_1$ to $M_2$,
and the sequence $\langle z_\beta :\,\beta \leq \alpha\,\rangle$
which is an increasing continuous sequence in $\niceT$.\\
Stage $\alpha$ in the game is as follows : \\
First, AIS chooses $z_\alpha$ of level $\alpha$ of $\niceT$, such that for every \mbox{$\beta <
\alpha \quad z_\alpha >^\niceT z_\beta$}. Then,
\begin{enumerate}
\item if $\alpha = 0$ then $f_\alpha = \emptyset$
\item if $\alpha$ is limit then $f_\alpha = \cup_{\beta <
\alpha}\,f_\beta$
\item if $\alpha = \beta + 1$ then  AIS chooses $A_1 \subseteq
M_1,\,A_2\subseteq M_2$ such that \mbox{$|A_1 \cup A_2| < 1+\mu$}.
Then ISO should choose $f_\alpha$ such that:\\
$f_\alpha$ is a partial isomorphism from $M_1$ to $M_2$,\hspace{1mm}$f_\beta
\subseteq f_\alpha$\\ 
\mbox{$A_1 \subseteq \Dom(f_\alpha),\,A_2 \subseteq \Rang(f_\alpha)$}
\end{enumerate}
The first player who cannot find a legal move loses the game. If
ISO has a winning strategy for  $\Game_{\niceT,\mu}(M_1,M_2)$, we
say that $M_1,M_2$ are \mbox{$EF_{\niceT,\mu}$ equivalent.}
\end{definition}
\begin{definition} \label{Sec1-Def-Square-F-lambda} We say that
$\boxtimes_{\niceF,\lambda}$ holds, if :
\begin{enumerate}
\item $\niceF$ is a set of partial functions from $\lambda$ to $\lambda$
\item if $f$ is a partial function from $\lambda$ to $\lambda$
then $f \in \niceF$ iff \\ for every countable \mbox{$u \subseteq
\Dom(f)$}\hspace{1mm}\quad$f \rest u \in \niceF$
\item there is no $f \in \niceF$ such that $\Dom(f) = \lambda$
\end{enumerate}
\end{definition}
\begin{definition}
If $\boxtimes_{\niceF,\lambda}$ holds, we define a tree
$\niceT_\niceF$ in the following way:
\begin{itemize}
\item the nodes are functions $f$ such that $f \in \niceF$ and
$\Dom(f)$ is an ordinal.
\item the order is inclusion
\end{itemize}
Note that this tree does not have a branch of length $\geq\lambda$
\end{definition}

\begin{remark}
\label{Sec1-Remark-Strength-Of-Game} If $\mathcal T$ is a tree of
size $\lambda$ with no $\lambda$-branch, we can assume without
loss of generality that $\niceT \subseteq \lambda$. Define
$\niceF$ by $f \in \niceF$ if $f$ is a partial function from
$\lambda$ to $\lambda$ such that $x < y \Rightarrow f(x) <^\niceT
f(y)$. We get that $\boxtimes_{\niceF,\lambda}$ holds, and
$\niceT$ can be embedded (as a partial order) in $\niceT_\niceF$
\end{remark}

\begin{theorem}
\label{Sec1-Theorem-Eq-T-lambda} Suppose :
\begin{enumerate}
\item $\cf(\lambda)=\lambda=\lambda^{\aleph_{0}}$
\item $\boxtimes_{\niceF,\lambda}$ holds
\item $\niceT = \niceT_\niceF$
\end{enumerate}
then :
\newline
There are non-isomorphic models $M_{1},M_{2}$ of size $\lambda$ which are \\
$EF_{\niceT,\lambda}$ equivalent.
\end{theorem}

\par \noindent Proof:\
First, we shall define a tool for constructing models.
\begin{definition}
${\mathfrak x}$ is a structure parameter if it consists of the
following objects:
\begin{enumerate}
\item a set $I$
\item a set $J_{s}$ for each $s \in I$, \st if $s_{1}\neq s_{2}$ then
$J_{s_{1}} \cap J_{ s_{2}} = \emptyset $.\\
denote $J=\bigcup_{s\in I}J_{s}$
\item sets $S,T$ \st $S\subseteq I\times I,\ T \subseteq J\times J$
\end{enumerate}
\end{definition}
\begin{definition}
For a given structure parameter ${\mathfrak x}$ we define a model
\mbox{$M=M_{{\mathfrak x}}$} in the following way :
\newline
First for each $s \in I$ let $\GG_{s}$ be an abelian group
generated freely by \mbox{$ \{ x_{t}:\,t\in J_{s} \}$} except of
the relation $\forall x(2x = 0)$. (We could have also used free
group or free abelian group, But this choice makes the proof a bit
simpler). We demand also that if $s_1 \neq s_2$ then
$\GG_{s_1}\cap \GG_{s_2} = \emptyset $.\\
For $(s_1 , s_2) \in S $, let $\GG_{s_1 , s_2} $ be the subgroup
of $\GG_{s_1}\times \GG_{s_2}$ generated by \mbox{$\{
(x_{t_1},x_{t_2}):\, (t_1 , t_2) \in T \cap ( J_{s_1} \times J_{s_2} )
\}$}.
\\
The universe of $M$ is $\bigcup_{s \in I}\,\GG_{s}$.
The vocabulary of $M$ consists of :
\begin{enumerate}
\item For each $a \in M$, \hspace{1pt} a unary function symbol $F_a$
\item For each $s \in I$, \hspace{1pt} a unary relation symbol $P_s$
\item For each $(s_1,s_2) \in S$, \hspace{1pt} a binary relation symbol
$Q_{s_1,s_2}$
\end{enumerate}
The interpretation of the symbols in $M$ is as follows :
\begin{enumerate}
\item For each $b \in M,\ s \in I, \ a \in \GG_{s}$ \hspace{1pt} if $b \in \GG_{s}$ then $F_a^M(b) =
a+b$ ; else $F_a^M(b) = b$
\item For each $s \in I$, $P_s^M = \GG_s$
\item For each $(s_1,s_2) \in S$,
$Q_{s_1,s_2}^M = \GG_{s_1,s_2} $
\end{enumerate}
\end{definition}
\begin{lemma}
\label{Sec1-LemmaCondOfBeingAuto}
 Suppose $I' \subseteq I$ and $f$ is a function,
$f:\,\bigcup_{s \in I'}\,\GG_s \rightarrow M$. Then $f$ is a
partial automorphism of $M$ iff :
\begin{enumerate}
\item for each $s \in I' \ \ f(0_{\GG_s}) \in \GG_s$
\item for each $s \in I',\ a \in \GG_s$  we have $f(a)=f(0_{\GG_s})+ a $
\item for each $s_1,s_2 \in I'$ if $(s_1,s_2) \in S $ then
$(f(0_{\GG_{s_1}}),f(0_{\GG_{s_2}})) \in \GG_{s_1,s_2} $
\end{enumerate}
\end{lemma}
\par \noindent
Proof:\ Suppose $f$ is a partial automorphism then :
\begin{enumerate}
\item for each $s \in I' \ \ 0_{\GG_s} \in \GG_s = P_s^M \Rightarrow f(0_{\GG_s}) \in P_s^M = \GG_s$
\item for each $s \in I',\ a \in \GG_s$ \hspace{1pt} $f(a) = f(F^M_a(0_{\GG_s}))= F^M_a(f(0_{\GG_s})) = f(0_{\GG_s})+ a$
\item for each $s_1,s_2 \in I'$ if $(s_1,s_2) \in S $ then
$(0_{\GG_{s_1}},0_{\GG_{s_2}}) \in \GG_{s_1,s_2}$ (because it's a
subgroup of $\GG_{s_1} \times \GG_{s_2}$) but $\GG_{s_1,s_2} =
Q_{s_1,s_2}^M$, therefore we have
$(f(0_{\GG_{s_1}}),f(0_{\GG_{s_2}})) \in \GG_{s_1,s_2} $
\end{enumerate}
Similar arguments show the other direction.
\hfill$\square_{\ref{Sec1-LemmaCondOfBeingAuto}}$
\\\par \noindent
Now we shall define a structure parameter ${\mathfrak x}$. Then
define $M = M_{\mathfrak x}$. Then we will choose elements
$a_*,b_* \in M$, define $M_1 = (M,a_*),\ M_2 = (M,b_*)$ and show
that $M_1,\ M_2$ are as required in theorem
\ref{Sec1-Theorem-Eq-T-lambda}.
\par \noindent
Let ${\mathfrak x} = {\mathfrak x}_{\lambda,\niceF}$ be the following structure parameter :
\begin{enumerate}
\item $I = [\lambda]^{\aleph_0}$
\item For $u \in I$, $J_u$ consists of the quadruples $t
= (u,g,h,\zeta) $ where :
\begin{enumerate}
\item $g,h$ are functions from $u$ into $\lambda$
\item $\zeta$ is a function from ${\rm supRang}(g)\cap u$ into
$\lambda$
\item $\zeta \in \niceF$
\item $g,h$ are weakly increasing
\item $g(x) =g(y) \Rightarrow h(x) = h(y)$
\item $h(x) > x $
\end{enumerate}
For $t = (u,g,h,\zeta)$ we will denote $u = u^t,\ g = g^t,\ h =
h^t, \ \zeta = \zeta^t$
\item $S = \{(u_1,u_2) : \ u_1,u_2 \in I \ {\mathrm{and}} \ u_1 \subseteq u_2 \}$
\item $T = \{(t_1,t_2) : \ t_1,t_2 \in J,\ u^{t_1} \subseteq
u^{t_2} ,\ g^{t_1} \subseteq g^{t_2} ,\ h^{t_1} \subseteq
h^{t_2},\ \zeta^{t_1} \subseteq \zeta^{t_2} \ \} $
\end{enumerate}
Let $M = M_{\lambda,\niceF} = M_{\mathfrak x}$ be the corresponding model. Note that
\mbox{$|I| = \lambda^{\aleph_0} =\lambda$} and for each $u \in I$,
\hspace{1pt} \mbox{$|J_u| = \lambda^{\aleph_0} =\lambda $},
therefore $||M|| = \lambda $. Define : \mbox{$a_* =
0_{\GG_\emptyset},\ b_* =
x_{(\emptyset,\emptyset,\emptyset,\emptyset)} $}. $M_1 = (M,a_*),\
M_2 = (M,b_*) $.
\begin{claim}
\label{Sec1-Claim-Equiv-Models} $M_1,M_2$ are
$EF_{\niceT,\lambda}$ equivalent
\end{claim}
\par \noindent
Proof:
\begin{definition}
\label{Sec1-Def-G-of-lambda} We define a set of functions
${\mathcal G} = {\mathcal G}(\lambda)$ with a partial order
$\leq^{{\mathcal G}}$ in the following way :
\begin{enumerate}
\item For an ordinal $\alpha < \lambda \ \ {\mat G}_\alpha $ is
the set of functions $g$  which satisfy :
\begin{enumerate}
\item $ g : \gamma \rightarrow \alpha $ , $\gamma < \lambda $
\item $g$ is weakly increasing
\end{enumerate}
\item ${\mat G} = \bigcup_{\alpha < \lambda}\,{\mat G}_\alpha $
\item For each $ g \in \niceG $ such that $\Dom(g)= \gamma$  we define $h_g
:\, \gamma \rightarrow \gamma+1 $ by : $h_g(x) = \Min(\{y : y <
\gamma \wedge g(y) > g(x)\}\cup \{\gamma\})$
\item $g_1 \leq^\niceG g_2 $ if \mbox{$g_1 \subseteq g_2$ and $h_{g_1}
\subseteq h_{g_2} $}
\end{enumerate}
\end{definition}
\begin{claim}
\label{Sec1-Claim-G-of-lambda}
\begin{enumerate}
\item $g(x) = g(y) \Rightarrow h_g(x) = h_g(y)$
\item $h_g(x) > x$
\item $h_g$ is weakly increasing
\item For every $g_1,g_2 \in  \niceG \ \ g_1 \leq^\niceG g_2 $ iff
\begin{enumerate}
\item $ \Dom(g_1) = \gamma_1 \leq \gamma_2 = \Dom(g_2) $ , $ g_1
\subseteq g_2$
\item if $\gamma_1 < \gamma_2$ then $g_2(\gamma_1) > g_2(x)$ for
every $x < \gamma_1$
\end{enumerate}
\item If $g_1 \in \niceG_\alpha $ and $\Dom(g_1) < \gamma <
\lambda $ then there is $g_2 \in \niceG_{\alpha+1}$ such that $g_1
\leq^\niceG g_2$ and $\Dom(g_2) = \gamma $
\item If $\delta < \lambda$ and we have $\langle g_\alpha
:\,\alpha < \delta \rangle $ such that $g_\alpha \in \niceG_\alpha
$ and \mbox{$\beta < \alpha \Rightarrow g_\beta \leq^\niceG
g_\alpha $,} then $g = \bigcup_{\alpha < \delta}\,g_\alpha$
satisfies $g \in \niceG_\delta$ and $g_\alpha \leq^\niceG g$ for
each $\alpha < \delta$
\end{enumerate}
\end{claim}
\par \noindent
Proof:\
\begin{enumerate}
\item[(1)-(3)] Easy
\item[(4)] If there is $x < \gamma_1$ such that $g_2(\gamma_1) = g_2(x)$
then \mbox{$h_{g_2}(x) = h_{g_2}(\gamma_1)> \gamma_1 \geq
h_{g_1}(x)$} so $g_1 \nless^\niceG g_2$. On the other direction,
if $g_1 \subset g_2$ and \mbox{$g_2(\gamma_1) > g_2(x)$} for every \mbox{$x <
\gamma_1$}, then for every such $x$: If there is $y < \gamma_1$ such
that $g_1(y) > g_1(x)$, let $y'$ be the minimal $y$ which
satisfies this. We get $h_{g_1}(x) = h_{g_2}(x) = y'$. If there is
no such $y$, we get \mbox{$h_{g_1}(x) = h_{g_2}(x) = \gamma_1$.}
Therefore we have $h_{g_1}\subset h_{g_2}$.
\item[(5)] Define $g_2 :\,\gamma\rightarrow \alpha+1$ by : \\
For $x\in \Dom(g_1) \ \ \ g_2(x) = g_1(x)$.\\
For $x \in \gamma\setminus \Dom(g_1) \ \ \ g_2(x) = \alpha$. \\
By (4) we get that $g_1 \leq^\niceG g_2$
\item[(6)] Remember that $\lambda$ is regular therefore $\bigcup_{\alpha < \delta}\,\Dom(g_\alpha) < \lambda $
\hfill$\square_{\ref{Sec1-Claim-G-of-lambda}}$\\
\end{enumerate}

\par \noindent
Now we will describe a winning strategy for ISO in the game
$\Game_{\niceT,\lambda}(M_1,M_2)$.\\
In stage $\alpha$ of the game ISO will choose a function
$g_\alpha$ such that :
\begin{enumerate}
\item $g_\alpha \in \niceG_\alpha$
\item $\beta < \alpha \Rightarrow g_\beta \leq^\niceG g_\alpha$
\item If $\alpha$ is a successor ordinal and in stage $\alpha$ AIS
chose the sets $A_1,A_2$ then for each $u \in I $ such that $(A_1
\cup A_2)\cap \GG_u \neq \emptyset$ we have \mbox{$u \subseteq
\Dom(g_\alpha)$}
\end{enumerate}
The choice of $g_\alpha$ is done in the following way :
\begin{enumerate}
\item $g_0 = \emptyset$
\item If $\alpha$ is limit, then $g_\alpha = \cup_{\beta <
\alpha}\,g_\beta$.\\ By \ref{Sec1-Claim-G-of-lambda} $g_\alpha \in
\niceG_\alpha$ and $\beta < \alpha \Rightarrow g_\beta  \leq
^\niceG g_\alpha$
\item If $\alpha = \beta + 1$ and in stage $\alpha$ AIS chose the sets
$A_1,A_2$, ISO will choose $\gamma < \lambda $ such that
$\Dom(g_\beta) < \gamma $ and $u \subseteq \gamma $ for every $u
\in I $ such that $(A_1 \cup A_2)\cap u \neq \emptyset $ (Such
$\gamma$ exists because $|A_1 \cup A_2| + \aleph_0 < \lambda )$ .
By \ref{Sec1-Claim-G-of-lambda} there is $g \in \niceG_\alpha$
such that $\Dom(g) = \gamma$ and $g_\beta \leq^\niceG g$. ISO will
choose such a function as $g_\alpha$ .
\end{enumerate}
Now remember that if $\alpha = \beta + 1$, then in stage $\alpha$
AIS has to choose a node on level $\alpha$, which is actually a
function $\zeta_\alpha: \, \alpha \rightarrow \lambda,\,\,
\zeta_\alpha \in \niceF $. Then he chooses $A_1\subset M_1, \, A_2
\subset M_2$ . Then ISO has to choose partial isomorphism
$f_\alpha$ from $M_1$ to $M_2$ such that $f_\beta \subseteq
f_\alpha, \, A_1 \subseteq \Dom(f_\alpha), \, A_2 \subseteq
\Rang(f_\alpha) $ (See \ref{Sec1-Def-Game-With-Tree}). So, ISO
chooses $g_\alpha$, and then defines $f_\alpha$ according to 
$f_\beta, A_1, A_2, g_\alpha, \zeta_\alpha$ in the following way :\\
\mbox{$\Dom(f_\alpha) = \Dom(f_\beta) \cup \bigcup\{ \GG_u :\, u
\in I, \,(A_1 \cup A_2) \cap \GG_u \neq \emptyset\}$}.\\ For each
$u \in I $ we have \mbox{$\GG_u \subseteq \Dom(f_\alpha)$} or
$\GG_u \cap \Dom(f_\alpha) = \emptyset$. \\ If \mbox{$\GG_u
\subseteq \Dom(f_\alpha)$} we define $f_\alpha(0_{\GG_u})= x_t $,
where \mbox{$t = (u,g_\alpha \rest u,h_{g_\alpha} \rest
u,\zeta_\alpha \rest (u \cap{\rm supRang}(g_\alpha \rest u)))$}\\
(Note that because $g_\alpha \in \niceG_\alpha$, we have
$\Rang(g_\alpha) \subseteq \alpha = \Dom(\zeta_\alpha)$).\\
For every $a \in \GG_u$ we define \mbox{$f_\alpha(a) =
f_\alpha(0_{\GG_u}) + a$}. By the construction we get that if
$(u_1,u_2) \in S$ then
$(f_\alpha(0_{\GG_{u_1}}),f_\alpha(0_{\GG_{u_2}})) \in
\GG_{u_1,u_2}$ (because the corresponding couple of $t$-ies lays
in $T$). Therefore by \ref{Sec1-LemmaCondOfBeingAuto} $f_\alpha$
is a partial automorphism of $M$. We also have :
\begin{enumerate}
\item For $\beta < \alpha$ $g_\beta \subseteq g_\alpha,\,
h_{g_\beta} \subseteq h_{g_\alpha},\, \zeta_\beta \subseteq
\zeta_\alpha$. Therefore $f_\beta \subseteq f_\alpha$.
\item For each $\alpha > 0$
\mbox{$f_\alpha(a_*) = f_\alpha(0_{\GG_\emptyset}) =
x_{(\emptyset,\emptyset,\emptyset,\emptyset)} = b_*$}. Therefore
$f_\alpha$ is a partial isomorphism from $M_1 = (M,a_*)$ into $M_2
= (M,b_*) $
\end{enumerate}
\hfill$\square_{\ref{Sec1-Claim-Equiv-Models}}$
\begin{claim}
\label{Sec1-ClaimNotIso} $M_1,M_2$ are not isomorphic.
\end{claim}
\par \noindent
Proof:\ It is enough to show that $M$ is rigid( = doesn't have
a non-trivial automorphism).\\
Assume toward contradiction that $f \neq id$ is an automorphism of
$M$. For each $ u \in I$ we define $c_u = f(0_{\GG_u})$. By
\ref{Sec1-LemmaCondOfBeingAuto}, for each $u \subseteq w \in I$ we
have $(c_u,c_w) \in \GG_{u,w}$. \\
For each $u \subset w \in I$ and $t = (w,g,h,\zeta) \in J_w$ we
define $\pi_{w,u}(t) \in J_u$ by \mbox{$\pi_{w,u}(t)=: (u,g \rest
u,h \rest u,\zeta \rest {\rm supRang}(g \rest u) \cap u )$}. By
the definition of $T$ we have that if $ t \in J_w , r \in J_u$
then $(r,t) \in T$ iff \mbox{$r = \pi_{w,u}(t)$}. We define
homomorphism $\hat{\pi}_{w,u}: \, \GG_w \rightarrow \GG_u$ by
$\hat{\pi}_{w,u}(x_t) = x_r$ where $r = \pi_{w,u}(t)$. We get that
$\GG_{u,w}$ is the subgroup of $ \GG_u \times \GG_w$ generated by
$\{(\hat{\pi}_{w,u}(x_t),x_t): \, t \in J_w \}$. Since $\{x_t: \,
t \in J_w \}$ generate $\GG_w$, we get that $\GG_{u,w} = \{
(\hat{\pi}_{w,u}(c),c): \, c \in \GG_w \}$.
\par \noindent
Define $n(u)$ to be the length of the reduced
representation of $c_u$ as a sum of the generators $\{x_t: \, t
\in J_u \}$. For $u \subseteq w \in I$ we get $n(u) \leq n(w)$
since $c_u = \hat{\pi}_{w,u}(c_w)$ and $\hat{\pi}_{w,u}$ sends one
generator to one generator. If for every $u \in I$ there is $w \in
I$ such that $n(w) > n(u)$ we can find a sequence $\langle u_n:\,
n<\omega \rangle$ such that $u_n \in I$ and $n(u_n) < n(u_{n+1})$.
Define $w = \cup_{n<\omega}\, u_n$, we get that $n(w)$ is infinite
- contradiction. Therefore, there is $u_* \in I $ such that
$n(u_*)$ is maximal. Since we assumed $f \neq id$ , $n(u_*) > 0$.\\
\par\noindent
Choose $t_* \in J_{u_*}$ such that $x_{t_*}$ appears in the
reduced representation of $c_{u_*}$. For each $u_* \subseteq w \in
I $ there is a unique $t(w) \in J_w$ such that $\pi_{w,u_*}(t(w))
= t_*$ and $x_{t(w)}$ appears in the reduced representation of
$c_w$. Such $t(w)$ exists because $c_{u_*} =
\hat{\pi}_{w,u_*}(c_w)$. It is unique because if there were two
such $t$-ies, $t_1,t_2$ then \mbox{$\hat{\pi}_{w,u_*}(x_{t_1}) =
\hat{\pi}_{w,u_*}(x_{t_2}) = x_{t_*}$}. Since in $\GG_{u_*} \
\forall x(2x=0)$ it implies $n(w) > n(u_*)$ which contradicts the
maximality of $n(u_*)$.\\
\par\noindent
Note that if $u\subseteq w\subseteq z\in I$ then $\pi_{z,u} =
\pi_{w,u}\circ \pi_{z,w} $. Therefore, by uniqueness of $t(w)$ if
$u_* \subseteq w \subseteq z \in I $ we have $t(w) =
\pi_{z,w}(t(z))$. For each \mbox{$u_* \subseteq w \in I $}, define \mbox{$
g^w = g^{t(w)}, \ h^w = h^{t(w)}, \ \zeta^w = \zeta^{t(w)}$}. If
$u_* \subseteq  w_1,w_2 \in I$ then the functions
$g^{w_1},h^{w_1},\zeta^{w_1}$ and $g^{w_2},h^{w_2},\zeta^{w_2}$
are respectively compatible, since $t(w_1) = \pi_{z,w_1}(t(z))$
and $t(w_2) = \pi_{z,w_2}(t(z))$ where $z = w_1 \cup w_2$. Define
$g = \cup \{ g^w :\, u_* \subseteq w \in I \} \\ h = \cup \{ h^w
:\, u_* \subseteq w \in I \} \\ \zeta = \cup \{ \zeta^w :\, u_*
\subseteq w \in I \}$.\\ We get:
\begin{enumerate}
\item $\Dom(g) = \Dom(h) = \lambda$
\item $g,h$ are weakly increasing
\item $h(x) > x$
\item $g(x) = g(y) \Rightarrow h(x) = h(y)$
\item $\zeta \in \niceF$ (this is by
\ref{Sec1-Def-Square-F-lambda}(2) )
\item $\supRang(g) \subseteq \Dom(\zeta)$
\end{enumerate}
By \ref{Sec1-Def-Square-F-lambda}(3) $\Dom(\zeta) \neq \lambda$.
Therefore by (6) \mbox{${\rm supRang}(g) < \lambda$}. Since $g$ is
weakly increasing and $\lambda$ is regular, there is $\alpha_0 <
\lambda $ such that for every \mbox{$ \alpha_0 < \alpha < \lambda
\ \ g(\alpha) = g(\alpha_0)$}. By (4) we get that for every
\mbox{$ \alpha_0 < \alpha < \lambda \ \ h(\alpha) = h(\alpha_0)$}.
Choose $\alpha > h(\alpha_0) > \alpha_0$ and get that $h(\alpha) <
\alpha$ contradicting(3).\hfill $\square_{\ref{Sec1-ClaimNotIso}}
\square_{\ref{Sec1-Theorem-Eq-T-lambda}}$
\newpage

%
%
\section{Games with trees for singular $\lambda =
\lambda^{\aleph_0}$} 
\noindent
It is clear that for $\lambda$ singular we
cannot expect the same result as in the previous section, since
the AIS player would be able to list all the members of $M_1,M_2$.
Thus, we prove a weaker result - we allow AIS to choose only one
element in each turn. We also remark in \ref{Section2-Generalize-866} that this result generalizes the 
result in \cite{SH-836} for such $\lambda$.

\begin{theorem}
\label{Sec2-Theorem-Eq-T,1} Suppose :
\begin{enumerate}
\item $\cf(\lambda) < \lambda=\lambda^{\aleph_{0}}$
\item $\boxtimes_{\niceF,\lambda}$ holds
\item $\niceT = \niceT_\niceF$
\end{enumerate}
then :
\newline
There are non-isomorphic models $M_{1},M_{2}$ of size $\lambda$
which are\\
$EF_{\niceT,1}$ equivalent.
\end{theorem}
\begin{remark}
\label{Section2-Generalize-866}
We can show that Theorem \ref{Sec2-Theorem-Eq-T,1} generalizes the result in \cite{SH-836} by choosing appropriate $\niceF$. The result there shows the existence of two non-isomorphic models of size $\lambda$ which are equivalent under every EF game of length $< \cf(\lambda)$, which consists of sub-games of 
length $ < \lambda$, such that AIS chooses the length of each sub-game before it starts, and in every sub-game he chooses one element in each move - see the definitions there. Now, an appropriate $\niceF$ can be chosen by looking at the proof there, but we will take a shortcut - we will use the result instead of the proof. Let us choose a pair of models $M_1,M_2$ as in the result in \cite{SH-836}. Without loss of generality assume that the universe of $M_1$ is $\lambda \times \{1\}$, and the universe of $M_2$ is $\lambda \times \{2\}$.
We can take $\niceF$ to be the set of functions $f$ which satisfy the following conditions:
\begin{enumerate}
\item $\Dom(f) \subseteq \lambda,\,\,\, \Rang(f) \subseteq \lambda $
\item define a partial function $f'$ from $M_1$ to $M_2$ by:
\begin{enumerate}
\item $\Dom(f') = \Dom(f) \times \{1\}$
\item for every $\alpha \in \Dom(f)$, $f'((\alpha,1)) = (f(\alpha),2)$ 
\end{enumerate}
then, $f'$ is a partial isomorphism.
\end{enumerate}
Now, it is not hard to see that  $EF_{\niceT_{\niceF},1}$ equivalence implies equivalence as in the result of \cite{SH-836}.
\end{remark}
\par \noindent
Proof of theorem \ref{Sec2-Theorem-Eq-T,1}: \\
Denote $\kappa = \cf(\lambda)$. ($\kappa > \aleph_0$ because
$\lambda = \lambda^{\aleph_0})$. Let $\langle \mu_i \, : i <
\kappa \rangle$ be an increasing and continuous sequence such that: 
$\mu_0 = 0$, \mbox{${\mu_i}^+ < \mu_{i+1} = \cf(\mu_{i+1})$}, $i > 0 \Rightarrow \mu_i > \aleph_0$,
 $\cup_{i<\kappa} \, \mu_i = \lambda $. For every
$\alpha < \lambda$ there is a unique $i < \kappa$, such that
$\alpha \in [\mu_i,\mu_{i+1})$. We denote $i = {\mathbf
i}(\alpha)$.
\par
We define a structure parameter ${\mathfrak x} = {\mathfrak x}_{\niceF,\lambda}$ in the following
way:
\begin{enumerate}
\item $I = [\lambda]^{\aleph_0}$
\item for $u \in I$ $J_u$ is the collection of quadruples
$t =(u,g,h,\zeta)$ such that :
\begin{enumerate}
\item $g,h$ are functions from $u$ into $\lambda$, $\zeta$ is a
function from some subset of $u$ into $\lambda$.
\item $\zeta \in \niceF$
\item for every $x \in u,\quad g(x) \in
[\mu_{\ibold(x)},\mu_{\ibold(x)}^+],\, h(x) \in
[\mu_{\ibold(x}),\mu_{\ibold(x)+1}]$

\item $g,h$ are weakly increasing
\item $g(x) = g(y) \Rightarrow h(x) = h(y)$
\item $h(x) > x$
\item \mbox{$\Dom(\zeta) = u \cap \bigcup\, \{\mu_{\ibold(x)} : x \in u \ {\rm and} \ h(x) = \mu_{\ibold(x)+1} \}$}
\end{enumerate}
For $t =(u,g,h,\zeta)$ we denote \mbox{$ u = u^t, \, g = g^t, \, h
= h^t, \, \zeta = \zeta^t$}
\item $S = \{ (u_1,u_2) : \, u_1,u_2 \in I, u_1 \subseteq u_2 \}$
\item \mbox{$T = \{(t_1,t_2) \in J : u^{t_1} \subseteq u^{t_2}, \, g^{t_1}
\subseteq g^{t_2}, \, h^{t_1} \subseteq h^{t_2}, \, \zeta^{t_1}
\subseteq \zeta^{t_2} \}$}
\end{enumerate}
Let $M = M_{\niceF,\lambda} = M_{\mathfrak x}$ be the corresponding model. Define :
\mbox{$a_* = 0_{\GG_\emptyset},\ b_* =
x_{(\emptyset,\emptyset,\emptyset,\emptyset)} $}. $M_1 = (M,a_*),\
M_2 = (M,b_*) $.
\begin{claim} \label{Sec2-Claim-Equiv-Models}
$M_1,M_2$ are $EF_{\niceT,1}$ equivalent
\end{claim}
\par \noindent
Proof:\\
\begin{definition}
\label{Sec2-Def-W} We define a partially ordered set of functions
$(\niceW,\leq^\niceW)$, which depends on the sequence $\langle
\mu_i \, : i < \kappa \rangle$ in the following way :
\begin{enumerate}
\item we define a set $\niceB$ such that $\bar{\beta} \in \niceB$
iff:
\begin{enumerate}
\item $\bar{\beta} = \langle \beta_i : \, i < \kappa \rangle$,
$\mu_i \leq \beta_i \leq \mu_{i+1}$
\item there is $j = \jbold(\bar{\beta}) < \kappa$ such that \mbox{$i <
\jbold(\bar{\beta})\Leftrightarrow \beta_i = \mu_{i+1}$}
\end{enumerate}
\item for $\bar{\beta} \in \niceB$ we define $\niceW_{\bar{\beta}}$
to be the set of functions $g$ which satisfy :
\begin{enumerate}
\item $\Dom(g) = \cup_{i< \kappa}\,[\mu_i,\beta_i)$
\item $g$ is weakly increasing
\item  for every $i < \kappa,\, x \in [\mu_i,\beta_i)$ we
have $g(x) \in [\mu_i,\mu^+_i]$, and if \mbox{$g(x) = \mu^+_i$}
then $i < \jbold(\bar{\beta)}$
\end{enumerate}
\item for $j < \kappa$ we define $\niceW_j =
\cup\{ \niceW_{\bar{\beta}} : \, \jbold(\beta) \leq j \}$
\item for $g \in \niceW_{\bar{\beta}}$ we define a function $h_g$ in the
following way :\\
 $\Dom(h_g) = \Dom(g)$ and for $i < \kappa, \,x \in
 [\mu_i,\beta_i)$ we define \\
 \mbox{$h_g(x)= \Min( \{ y : \mu_i \leq y < \beta_i \wedge g(y) > g(x) \} \cup \{\beta_i\})$}
\end{enumerate}
\end{definition}
\begin{claim}
\label{Sec2-Claim-W}
\begin{enumerate}
\item $g(x) = g(y) \Rightarrow h_g(x) = h_g(y)$
\item $h_g(x) > x $
\item $h_g$ is weakly increasing
\item $ x \in [\mu_i,\mu_{i+1}) \Rightarrow h_g(x) \in
[\mu_i,\mu_{i+1}]$
\item Suppose that $g_1 \in \niceW_{\bar{\beta}^1}$ , $g_2 \in \niceW_{\bar{\beta}^2}$ then $g_1 \leq^\niceW g_2 $ iff
\begin{enumerate}
\item $ g_1 \subseteq g_2 $ (therefore for every $ i < \kappa$
$\beta^1_i \leq \beta^2_i$ )
\item for every $i < \kappa$\hspace{1mm} if $\beta^1_i < \beta^2_i$\\
then for every $x \in [\mu_i,\beta^1_i),\,\,g_2(x) <
g_2(\beta^1_i)$
\end{enumerate}
\item if $g_1 \in \niceW_j$ and $\bar{\beta} \in \niceB$,
$\jbold(\bar{\beta}) \leq j$ then there is $g_2 \in \niceW_j$ such
that $g_1 \leq^\niceW g_2 $ and $\cup_{i <\kappa} \,
[\mu_i,\beta_i) \subseteq \Dom(g_2)$
\item if  $\delta < \mu^+_j$ and $\langle g_\alpha \, : \alpha <
\delta \rangle $ satisfy $\alpha < \beta \Rightarrow g_\alpha \leq
^\niceW g_\beta$, $g_\alpha \in \niceW_j$ then there is $g \in
\niceW_j$ such that $\alpha < \delta \Rightarrow g_\alpha
\leq^\niceW g $
\end{enumerate}
\end{claim}
\par \noindent
Proof : \begin{enumerate}
\item[(1) - (4)] easy
\item[(5)] like in the proof of \ref{Sec1-Claim-G-of-lambda}
\item[(6)] We may assume that $\Dom(g_1) \subseteq \cup_{i<
\kappa}\,[\mu_i,\beta_i)$. Define for $i < \kappa$ \mbox{$\gamma_i
= \mu_i + \Sup \{g_1(x) : x \in \Dom(g_1) \cap
[\mu_i,\mu_{i+1})\}$}. Since $g_1 \in \niceW_j$ we have $i \geq j
\Rightarrow \gamma_i < \mu^+_i $. Define for $i < \kappa$
\begin{displaymath}
\gamma^*_i = \left\{ \begin{array}{ll} \mu^+_i \quad \textrm{if}
\quad i < j
\\
\gamma_i \quad \textrm{if} \quad i \geq j
\end{array} \right.
\end{displaymath}
Now define $g_2$ by : $\Dom(g_2) = \cup_{i<
\kappa}\,[\mu_i,\beta_i)$, and for every $i < \kappa$, \mbox{$x
\in [\mu_i,\beta_i)$} we define :
\begin{displaymath}
g_2(x) = \left\{ \begin{array}{ll} g_1(x) \quad \textrm{if} \quad
x \in \Dom(g_1)
\\
\gamma^*_i \quad \textrm{if} \quad x \notin \Dom(g_1)
\end{array} \right.
\end{displaymath}
Since $\jbold(\bar{\beta}) \leq j$ we have $g_2 \in \niceW_j$. By
(5) we have $g_1 \leq^\niceW g_2$.
\item[(7)] Define for every $i < \kappa$ :\\
$\beta_i = \sup(\,\cup_{\alpha < \delta}\,\Dom(g_\alpha) \cap
[\mu_i,\mu_{i+1})\,) + \mu_i$\\
$\gamma_i = \sup( \,\cup_{\alpha < \delta}\,\Rang(g_\alpha \rest
[\mu_i,\mu_{i+1}) ) + \mu_i$\\
For every $\alpha < \delta \quad g_\alpha \in \niceW_j$. Therefore
for every $i \geq j$
\begin{itemize}
\item $\sup(\Dom(g_\alpha)\cap[\mu_i,\mu_{i+1}) ) < \mu_{i + 1}$
\item $\supRang(g_\alpha \rest [\mu_i,\mu_{i+1}) ) < \mu^+_i$.
\end{itemize}
Therefore, since \mbox{$\delta < \mu^+_j \leq \mu^+_i < \mu_{i+1}
= \cf(\mu_{i+1})$}, we get that for \mbox{$ i \geq j$}
\mbox{$\beta_i <
\mu_{i+1}$ and $\gamma_i < \mu^+_i$}.\\
Define for $i <\kappa$ : $\beta^*_i = \left\{
\begin{array}{ll}
\mu_{i+1} \quad i < j \\
\beta_i  \quad i \geq j
\end{array} \right.$
$\gamma^*_i = \left\{
\begin{array}{ll}
\mu^+_i \quad i < j \\
\gamma_i + 1 \quad i \geq j
\end{array} \right.$\\
Denote $g' = \cup_{\alpha < \delta}\, g_\alpha$.\\
Define $g \in \niceW_j$ by :\\
$\Dom(g) = \cup_{i < \kappa}\,[\mu_i,\beta^*_i)$\\
For $ i < \kappa, \, x \in [\mu_i,\beta^*_i) \quad g(x) = \left\{
\begin{array}{ll}
g'(x) \quad x \in \Dom(g') \\
\gamma^*_i  \quad x \notin \Dom(g')
\end{array} \right. $\\
By (5) we get that $\alpha < \delta \Rightarrow g \geq^\niceW
g_\alpha$.\hfill $\square_{\ref{Sec2-Claim-W}}$\\
\end{enumerate}
\par \noindent
Now we will describe a winning strategy for ISO :\\
In every stage $\alpha$ in the game ISO will choose a function
$g_\alpha$ such that :
\begin{enumerate}
\item $g_\alpha \in \niceW_{\ibold(\alpha)+1}$
\item $\eps < \alpha \Rightarrow g_\eps \leq^\niceW g_\alpha$
\item If in stage $\alpha$ AIS chose an element from $\GG_u$
then $u \subseteq \Dom(g_\alpha)$
\end{enumerate}
ISO can choose such $g_\alpha$ in the following way :
\begin{enumerate}
\item for $\alpha = 0 \quad g_0 = \emptyset$
\item for $\alpha$ limit, since $\alpha < \mu_{\ibold(\alpha)+1}$
and for every $\eps < \alpha \quad g_\eps \in
\niceW_{\ibold(\alpha)+1}$, we can use (7) of \ref{Sec2-Claim-W}.
\item If $\alpha = \eps +1$ and in stage $\alpha$ AIS chose
element from $\GG_u$,then we choose $\bar{\beta} = \langle \beta_i
\, : i < \kappa \rangle$ in the following way : \\
If $i < \ibold(\alpha) + 1$ then $\beta_i = \mu_{i+1}$. Else,
$\mu_{i+1} > \alpha$. we choose \mbox{$\beta_i < \mu_{i+1}$} such
that $ u \cap [\mu_i,\mu_{i+1}) \subseteq [\mu_i,\beta_i)$. Now
$\jbold(\bar{\beta}) = \ibold(\alpha)+1$, so by (6) of
\ref{Sec2-Claim-W} we can find $g \in \niceW_{\ibold(\alpha)+1}$
such that $g_\eps \leq^\niceW g$ and \mbox{$\cup_{i < \kappa}\,
[\mu_i,\beta_i) \subseteq \Dom(g)$}. Define $g_\alpha = g$.
\end{enumerate}
Now if $\alpha = \eps+1$ and in stage $\alpha$ AIS chose an
element from $\GG_u$ and the node $\zeta_\alpha \in \niceT$, then
ISO will define the automorphism $f_\alpha$
according to $g_\alpha,\zeta_\alpha$ :\\
$\Dom(f_\alpha) = \Dom(f_\eps) \cup \GG_u$. For every $w$ such
that $\GG_w \subseteq \Dom(f_\alpha)$,\\
\mbox{$f_\alpha(0_{\GG_w}) = x_t$} where $t = (w,g_\alpha \rest
w,h_{g_\alpha} \rest w,\zeta_\alpha
\rest v)$\\
where $v = w \cap \{\mu_{\ibold(x)}:\,x \in w \wedge
h_{g_\alpha}(x) =\mu_{\ibold(x)+1}\}$\\
(Note that $v \subseteq \alpha=\Dom(\zeta_\alpha)$, because
$g_\alpha \in \niceW_{\ibold(\alpha) + 1}$ )\\
As in section 1, we get that $f_\alpha$ is a partial isomorphism
and $\eps < \alpha \Rightarrow f_\eps \subseteq f_\alpha$.
\hfill$\square_{\ref{Sec2-Claim-Equiv-Models}}$

\begin{claim}
\label{Sec2-ClaimNotIso} $M_1,M_2$ are not isomorphic.
\end{claim}
\par \noindent
Proof:\ We imitate the proof of \ref{Sec1-ClaimNotIso}. It is
enough to show that $M$ is rigid. Assume toward contradiction that
$f \neq id$ is an automorphism of $M$. For each $u \subset w \in
I$ and $t = (w,g,h,\zeta) \in J_w$ we define $\pi_{w,u}(t) \in
J_u$ by $\pi_{w,u}(t) =(u,g^t\rest u,h^t\rest u,\zeta^t\rest v)$ \\
where \mbox{$v = \cup\{\mu_{\ibold(x)}\, :x \in u \wedge h^t(x) =
\mu_{\ibold(x)+1}\, \}\cap u $}.\\
We proceed as in the proof of \ref{Sec1-ClaimNotIso}, and we get
that we can find functions $g,h,\zeta$ such that:
\begin{enumerate}
\item $\Dom(g) = \Dom(h) = \lambda, \,\, \Dom(\zeta) \subseteq \lambda$
\item if $\ibold(x) = i$ then $g(x) \in [\mu_i,\mu^+_i],\,\, h(x) \in [\mu_i,\mu_{i+1}]$
\item $g,h$ are weakly increasing
\item $g(x) = g(y) \Rightarrow h(x) = h(y)$
\item $h(x) > x$
\item $h(x) = \mu_{\ibold(x)+1} \Rightarrow \mu_{\ibold(x)}
\subseteq \Dom(\zeta)$
\item $\zeta \in \niceF$
\end{enumerate}
By (7) we get that $\Dom(\zeta) \neq \lambda$, therefore by (6)
there is $i < \kappa$ such that such that $\ibold(x) = i
\Rightarrow \ibold(h(x)) = i$. By (2) $\ibold(x) = i \Rightarrow
g(x) \leq \mu^+_i$. By (3) $g$ is weakly increasing. Since
$\mu_{i+1}= \cf(\mu_{i+1}) > \mu^+_i$, we can find $\alpha_0$ such
that $\alpha_0 \leq x < \mu_{i+1} \Rightarrow g(x) = g(\alpha_0)$.
By (5) $h(\alpha_0) > \alpha_0$. By the choice of $i$ we get that
$h(\alpha_0) < \mu_{i+1}$. Choose $h(\alpha_0) < x < \mu_{i+1}$.
We get $h(x) > x > h(\alpha_0)$ but $g(x) = g(\alpha_0)$. This
contradicts (4). Therefore we proved that $M$ is rigid.
\hfill$\square_{\ref{Sec2-ClaimNotIso}}\square_{\ref{Sec2-Theorem-Eq-T,1}}$
\newpage

%
%
\section{$\lambda$ regular $ > \beth_\omega$}
\noindent
In this section we show a result which holds for every $\lambda$
regular $> \beth_\omega$. In the previous sections we used the
assumption $\lambda = \lambda^{\aleph_0}$. Here we use instead of
it the existence of a set $\niceP \subset [\lambda]^{\aleph_0}$ of
size $\lambda$ which is "dense". By "dense" we mean that for every $A \in [\lambda]^{\beth_\omega}$
there is $B \subset A,\,B \in \niceP$.
\begin{remark}
\begin{enumerate}
\item Looking at the proof, one can see that instead of \mbox{$\lambda > \beth_\omega$},\\ 
it is enough to assume the following:
\begin{enumerate}
\item $\lambda > 2^{\aleph_0}$
\item There is $\niceP \subset [\lambda]^{\aleph_0}$ such that:
\begin{enumerate}
\item $|\niceP| = \lambda$
\item For every $A \in [\lambda]^\lambda$, there is $B \in \niceP$, $B \subset A$.
\end{enumerate}
\end{enumerate}
\item It is possible that it can be proved in ZFC that every $\lambda > 2^{\aleph_0}$ satisfies (1)(b) (It is a problem in cardinal arithmetic).
\end{enumerate}
\end{remark}

\begin{theorem}
\label{Sec3-Theorem-Equiv} Suppose:
\begin{enumerate}
\item $\lambda = \cf(\lambda) > \beth_\omega$
\item $\niceT$ is a tree of size $\lambda$  with no branch
of length $\lambda$
\end{enumerate}
Then : there are models $M_1,M_2$ of size $\lambda$ which are
$EF_{\niceT,\lambda}$ equivalent but not isomorphic
\end{theorem}
\par \noindent
Proof:
Let $\chi$ be large enough cardinal (for example $\chi =\beth_7(\lambda)$ ).
\begin{claim}
\label{Sec3-Claim-exists-inner-model}
We can find $\Mfrak$ such that :
\begin{enumerate}
\item $\Mfrak$ is elementary sub-model of $\niceH(\chi)$
\item $\lambda+1 \subseteq \Mfrak$
\item $||\Mfrak|| = \lambda$
\item for every $\langle\,(x_i,z_i) : i < \lambda \rangle$ such that $x_i \in \Mfrak$, $z_i \in \niceT$ for every 
      $i < \lambda$ there is an increasing sequence $\langle i_n : n < \omega \rangle$ such that:
      \begin{enumerate}
      \item $\langle (x_{i_n}, z_{i_n}) : n < \omega \rangle \in \Mfrak$
      \item if in addition, for $i < j < \lambda$ the level of $z_i$ (in $\niceT$) is strictly less then the 
      level of $z_j$, then $\langle z_{i_n} : n < \omega \rangle$ is an antichain in the order $\leq^\niceT$
      \end{enumerate}
\end{enumerate}
\end{claim}
\par \noindent
Proof:\\
We use part of the RGCH theorem(see Shelah \cite{SH-RGCH})
\begin{rgchtheorem}
\label{RGCH-Theorem}
if $\lambda \geq \beth_\omega$ then there is regular $\kappa < \beth_\omega$
 and $\niceP \subseteq [\lambda]^{<\beth_\omega} $ such that :
\begin{enumerate}
\item $|\niceP| = \lambda$
\item for every $A \in [\lambda]^{\beth_\omega}$, We can find $\langle A_i :\, i < \epsilon \rangle$ such that:\\
$\epsilon < \kappa,\,\, A_i \in \niceP$ for every $i < \epsilon$ and $A = \bigcup_{i < \epsilon}\,A_i$
\end{enumerate}
\end{rgchtheorem}
\begin{corollary}
\label{RGCH-Corollary}
if $\lambda \geq \beth_\omega$ then we can find a set $\niceP^* \subseteq [\lambda]^{\aleph_0}$ such that
$|\niceP^*| = \lambda$ and for every $A \in [\lambda]^{\beth_\omega}$ there is $B \in \niceP^*$ such that
$B \subseteq A$
\end{corollary}
\noindent
Proof:\\
Choose $\kappa$ and $\niceP$ as in \ref{RGCH-Theorem} and define $\niceP^* = \bigcup\{ [A]^{\aleph_0} :\, A \in \niceP\}$
\hfill$\square_{\ref{RGCH-Corollary}}$\\

\par \noindent
We construct $\Mfrak_n$ for every $n < \omega$ such that:
\begin{enumerate}
\item $\Mfrak_0$ is an elementary sub-model of
$\niceH(\chi)$ such that $||\Mfrak_0|| = \lambda$,\\
$\lambda+1 \subseteq \Mfrak_0$ and for every $A \in [\lambda]^{\beth_\omega}$, there is
$B \in \Mfrak_0 \cap [\lambda]^{\aleph_0},\,\ B
\subset A$ (This is possible by \ref{RGCH-Corollary} )
\item $||\Mfrak_n|| = \lambda$
\item $\Mfrak_n$ is an elementary sub-model of $\niceH(\chi)$
\item if $A \in \Mfrak_n$ and $|A| \leq \lambda$, then $A \subseteq \Mfrak_{n+1}$
\item $\Mfrak_n \in \Mfrak_{n+1}$, $\Mfrak_n \subset \Mfrak_{n+1}$
\end{enumerate}

\par \noindent 
Now, let $\Mfrak = \bigcup_{n < \omega}\,\Mfrak_n$. We will prove that $\Mfrak$ satisfies the conclusion of claim \ref{Sec3-Claim-exists-inner-model}.\\
Suppose that \mbox{$\langle\,(x_i,z_i) : i < \lambda \rangle \subseteq \Mfrak \times \niceT$} satisfies $x_i \in \Mfrak$, $z_i \in \niceT$ for every $i < \lambda$. 
We may assume without loss of generality, that there is $n_0 < \omega$ such that $\{(i,x_i,z_i) : i < \lambda \} \subseteq \Mfrak_{n_0}$.
If the condition in \ref{Sec3-Claim-exists-inner-model} 4(b) is not satisfied, then we are done, because we can 
find $A \in  [\lambda]^{\aleph_0}$ such that \mbox{$\{(i,x_i,z_i) : i \in A \} \in \Mfrak_{n_0+1}$}. 
(Because in $\Mfrak_{n_0 + 1}$ there is one to one correspondence between $\lambda \times \Mfrak_{n_0} \times \niceT$ and $\lambda$, and every subset of $\lambda$ of size $\beth_\omega$ has infinite countable subset that is a member of $\Mfrak_0$).\\
If the condition in \ref{Sec3-Claim-exists-inner-model} 4(b) is satisfied, then we have 2 cases:\\
case (1):\\
We can find $A \in [\lambda]^{\beth_\omega}$ such that $\langle z_i :\, i \in A \rangle$ is an antichain in $\leq^\niceT$\\
case (2):\\
We cannot find such $A$.\\
If we are in case(1) then we are done in the same way as before.\\
Suppose we are in case(2):
\begin{claim} \label{Sec3-small-1} for every $j < \lambda$, we can find $j < i_0 < i_1 < i_2 < \lambda$, such that $z_{i_0} <^\niceT z_{i_1},z_{i_2}$ and $z_{i_1},z_{i_2}$ are not comparable in $\leq^\niceT$.
\end{claim}
\noindent
Proof: assume toward contradiction that there is $j* < \lambda$, such that we can't find $j* < i_0 < i_1 < i_2 < \lambda$ which are as in the claim.\\ Define \mbox{$C = \{z_i \, :  j*< i <  \lambda\}$}. Then, being comparable in $\leq^\niceT$ is an equivalence relation on $C$. Since $\lambda$ is regular, either there are $\lambda$ equivalence classes or there is an equivalence class of size $\lambda$. In other words,\\ $C$ contains an antichain or a chain of size $\lambda$, both options are not possible, the first since we are in case (2) and the second since $\niceT$ doesn't have a $\lambda$ branch -contradiction.\hfill$\square_{\ref{Sec3-small-1}}$
\par \noindent
\\By claim \ref{Sec3-small-1} we can choose for every $ j < \lambda$ a triple $i_0(j),i_1(j),i_2(j)$
such that :
\begin{enumerate}
\item $i_0(j) < i_1(j) < i_2(j) < \lambda$
\item $ j < j' \Rightarrow i_2(j) < i_0(j')$
\item $z_{i_0(j)} <^\niceT z_{i_1(j)},z_{i_2(j)}$
\item $z_{i_1(j)}$ and  $z_{i_1(j)}$ are not comparable in $\leq^\niceT$
\end{enumerate}
We choose $A \in [\lambda]^{\aleph_0}$ , such that
 \mbox{$\{(j,i_0(j),i_1(j),i_2(j),x_j,z_j):\,j \in A\} \in \Mfrak_{n_0 +1}$}
Using Ramesy theorem in $\Mfrak_{n_0+1}$, we can find an increasing sequence \\
$\langle j_n:\, n < \omega \rangle$ such that:
\begin{enumerate}
\item for every $n < \omega, \,\,\, j_n \in A$
\item $\langle j_n:\, n < \omega \rangle \in \Mfrak_{n_0 +1}$
\item $\{ z_{i_1(j_n)} :\, n < \omega \} $ is a chain or an antichain in $\niceT$
\item $\{ z_{i_2(j_n)} :\, n < \omega \} $ is a chain or an antichain in $\niceT$
\end{enumerate}
Now we are done, since either $\{ z_{i_1(j_n)} :\, n < \omega \} $ or $\{ z_{i_1(j_n)} :\, n < \omega \} $ must be an antichain. Because if both are chains, we get that\\
\mbox{$z_{i_1(j_0)} <^\niceT z_{i_1(j_1)},\,\,z_{i_2(j_0)} <^\niceT z_{i_2(j_1)}$}. Since $z_{i_0(j_1)}$ is on higher level then $z_{i_1(j_0)},\,z_{i_2(j_0)}$ and it is $ <^\niceT z_{i_1(j_1)},\,z_{i_2(j_1)}$ We get that \mbox{$z_{i_1(j_0)},z_{i_2(j_0)} <^\niceT z_{i_0(j_1)}$} \\- contradiction, since by the construction they are not comparable.\hfill$\square_{\ref{Sec3-Claim-exists-inner-model}}$
\par \noindent \\
We choose $\Mfrak$ as in claim $\ref{Sec3-Claim-exists-inner-model}$.\\
We define a structure parameter $\xfrak = \xfrak(\Mfrak)$ in the
following way:
\begin{definition}
\label{Sec3-Def-Parameter}
\begin{enumerate}

\item $I$ consists of the objects of the form $(u,\Lambda)$ where:
\begin{enumerate}
\item $u \in \lambda^{<\aleph_0}$
\item $\Lambda \in \Mfrak$, $| \Lambda | \leq \aleph_0$, $\Lambda$ is a set of
partial functions with finite domain, from $\lambda$ to $\lambda$.
\end{enumerate}
for $s = (u,\Lambda)$ we denote $u = u^s,\, \Lambda = \Lambda^s$ .\\
We define $\Gamma(s) = u^s \cup \,\bigcup\{ \Dom(f):\, f \in
\Lambda^s\}$. Note that this a countable set.

\item For $s = (u,\Lambda) \in I$, $J_s$ consists of all the objects of the form \mbox{$t
= (u,\Lambda,g,h,F,z)$} where:
\begin{enumerate}
\item $g,h$ are functions from $u$ to $\lambda$
\item $F$ is a function from $\Lambda^2$ to $\{0,1\}$
\item $z \in \niceT$
\item Let $\alpha$ be the level of $z$ in the tree $\niceT$. then
$\alpha$ is minimal under the condition $\alpha > y $ for every
$y$ such that: \\
$ y \in \Rang(g)$ or there are $f_1,f_2 \in \Lambda$ such that
$F(f_1,f_2) =1$ and $y \in \Rang(f_1)$
\item There is a witness $(\gbold,\hbold)$ for $t$, which means
that :
\begin{enumerate}
\item $\Dom(\gbold) = \Dom(\hbold) \subseteq
\lambda,\,\,\Rang(\gbold) \cup \Rang(\hbold) \subseteq \lambda$
\item $\Gamma(s) \subseteq \Dom(\gbold)$
\item $\gbold,\hbold$ \hspace{1mm} are weakly increasing
\item $\hbold(x) > x$
\item $\gbold(x) = \gbold(y) \Rightarrow \hbold(x) = \hbold(y)$

\item $g \subseteq \gbold, h \subseteq \hbold$
\item for every $(f_1,f_2) \in \Lambda^2$\\
$F(f_1,f_2) = 1$ iff \hspace{1mm} $f_1 \subseteq \gbold \wedge\,\,
f_2 \subseteq \hbold$

\end{enumerate}
\end{enumerate}

\item $S = I^2$

\item $T$ consists of the pairs $(t_1,t_2) \in J^2$ where :
\begin{enumerate}
\item $t_1,t_2$ have a common witness
\item $z^{t_1},z^{t_2}$ are comparable in the order $\leq^\niceT$
\end{enumerate}

\end{enumerate}
\end{definition}
\begin{fact}
if:
\begin{enumerate}
\item $ s \in I, z \in\niceT$
\item $\gbold,\hbold$ satisfy conditions (i)-(v) from \ref{Sec3-Def-Parameter}(e)
\item $\Dom(\gbold) \subset \alpha$ where $\alpha$ is the level of $z$
\end{enumerate}
then:
\begin{enumerate}
\item there is unique $t \in J_s$ such that $(\gbold,\hbold)$ is a
witness for t, and $z^t \leq^\niceT z$. we denote $t =
t(s,\gbold,\hbold,z)$
\item if :
\begin{enumerate}
\item $\gbold',\hbold',z'$ also satisfy the conditions in (1)
\item $z,z'$ are comparable in $\leq^\niceT$
\item  $\gbold',\hbold'$ are compatible with $\gbold,\hbold$ respectively
\end{enumerate}
then :\\
$t(s,\gbold,\hbold,z) = t(s,\gbold',\hbold',z')$ \\
\end{enumerate}
\end{fact}
\noindent
Let $M = M_{\mathfrak x}$ be the corresponding model. We can check
that $||M|| = \lambda$. Let $a_* =
0_{\GG_{(\emptyset,\emptyset)}},\quad b_* =
x_{(\emptyset,\emptyset,\emptyset,\emptyset,\emptyset,z_*)}$ where $z_*$ is the
root of $\niceT$ (without loss of generality there is a root).\\
Define $M_1 = (M,a_*), M_2 = (M,b_*)$.

\begin{claim}
\label{Sec3-Claim-Equiv} $M_1,M_2$ are $EF_{\niceT,\lambda}$
equivalent.
\end{claim}
\par \noindent
We describe a winning strategy for ISO - this is very similar to
the proof of  \ref{Sec1-Claim-Equiv-Models}, so we will omit the
details. We are using the definitions in \ref{Sec1-Def-G-of-lambda} .\\
In every stage $\alpha$ of the game ISO will choose a function $\gbold_\alpha$ such that :
\begin{enumerate}
\item $\gbold_0 = \emptyset$
\item $\gbold_\alpha \in \niceG_\alpha$ (See definition of $\niceG_\alpha$ and $\leq^{\niceG}$ in  \ref{Sec1-Def-G-of-lambda} )
\item $\beta < \alpha \Rightarrow \gbold_\beta \leq^\niceG \gbold_\alpha$
\item If in stage $\alpha$ AIS chose the sets $A_1,A_2$ then\\
for each $s \in I$, if $\GG_s \cap (A_1 \cup A_2) \neq \emptyset$
then $\Gamma(s) \subseteq \Dom(\gbold_\alpha)$
\end{enumerate}
Now if $\alpha = \beta + 1$ and in stage $\alpha$ AIS chose the
sets $A_1,A_2$ and the node $z_\alpha$, ISO will define
$\hbold_\alpha = h_{\gbold_\alpha}$ and then define $f_\alpha$ by
:
\begin{enumerate}
\item $\Dom(f_\alpha) = \bigcup\{\GG_s : \Gamma(s) \subseteq
\Dom(\gbold_\alpha)\}$
\item for each $s$ such that $\GG_s
\subseteq \Dom(f_\alpha)$, $f_\alpha(0_{\GG_s}) = x_t$,\\
where $t =
t(s,\gbold_\alpha,\hbold_\alpha,z_\alpha)$\hfill$\square_{\ref{Sec3-Claim-Equiv}}$
\end{enumerate}

\begin{claim}
\label{Sec3-Claim-Not-Iso} $M_1,M_2$ are not isomorphic
\end{claim}
\par \noindent
Proof:\\
It is enough to show that $M$ is rigid. Assume toward
contradiction that $f \neq id $ is an automorphism of $M$. Denote
for $s \in I,\, c_s = f(0_{\GG_s})$. Denote $W_s = \{t \in J_s :\,
x_t\,\, \textrm{is in the reduced representation of}\,\, c_s \}$.
Since $f \neq id$ there is  $s^* = (u^*,\Lambda^*)$ such that
$W_{s^*} \neq \emptyset$. Note also that if $u^{s^*} \subseteq
u^s$ and $\Lambda^* \subseteq \Lambda^s$, then there is a natural
projection $\pi_{s,s^*}$ from $J_s$ into $J_{s^*}$ such that
$W_{s^*}\subseteq \Rang(\pi_{s,s^*}\rest\, W_s)$ (see the proof of
\ref{Sec1-ClaimNotIso}) therefore $W_s \neq \emptyset$.
\par
Choose for $i < \lambda \quad s_i,t_i,\alpha_i$ such that :
\begin{enumerate}
\item $s_i \in I,\,s_i = (u^* \cup \{\alpha_i\},\Lambda^*)$
\item $t_i \in W_{s_i}$
\item $\alpha_i < \lambda$
\item $i < j \Rightarrow h^{t_i}(\alpha_i) < \alpha_j$\\
\end{enumerate}
Case (*1) : $\Sup\{g^{t_i}(\alpha_i) :\, i < \lambda\} = \lambda$.
Then , since the level of $z^{t_i}$ in $\niceT$ must be greater then
$g^{t_i}(\alpha_i)$, we may assume that if $i < j$ then the level of $z^{t_i}$ is strictly less
then the level of $z^{t_j} $.\\
Case (*2): $\Sup\{g^{t_i}(\alpha_i) :\, i < \lambda\} < \lambda$.
Then by regularity of $\lambda$, we may assume that for every $i,j
<\lambda\quad g^{t_i}(\alpha_i) = g^{t_j}(\alpha_j)$\\
\par \noindent
Now, no matter in which case we are, we proceed in the following
way:\\
By the properties of $\Mfrak$ (see claim \ref{Sec3-Claim-exists-inner-model}) we can find a set $A \subset
\lambda$ such that:
\begin{enumerate}
\item $|A| = \aleph_0$
\item  $\{W_{s_i}:\, i \in A \} \in \Mfrak$
\item if we are in case (*1) $\{z^{t_i} :\, i \in A\}$ is an antichain (We can have that because in case(*2) the level
of $z^{t_i}$ is strictly increasing with $i$ - See \ref{Sec3-Claim-exists-inner-model} )
\end{enumerate}
We define $s^+ = (u^*,\bigcup_{i \in A}\,W_{s_i}\,\cup \Lambda^*)$. (Note that $\bigcup_{i \in A}\,W_{s_i} \in \Mfrak$,\\
therefore $s^+ \in I$)
\begin{claim}
\label{Sec3-Claim-unique-projection} For every $i \in
A$,\hspace{1mm} if : $r \in J_{s^+},\,t \in W_{s_i}, (r,t) \in
T$\\
then :
\begin{enumerate}
\item if $(\gbold,\hbold)$ is a witness for $r$ then $g^t
\subseteq \gbold,\,h^t \subseteq \hbold$
\item if $t \neq t' \in J_{s_i} $ then $(r,t') \notin T$
\end{enumerate}
\end{claim}
\par \noindent
Proof:
\begin{enumerate}
\item Let $(\gbold_0,\hbold_0)$ be a common witness for $r,t$. Then $g^t
\subseteq \gbold_0,\,h^t \subseteq \hbold_0$. Now $g^t,h^t \in
\Lambda^{s^+}$ therefore $(g^t,h^t) \in \Dom(F^r)$. since
$(\gbold_0,\hbold_0)$ is a witness for $r$ and $g^t \subseteq
\gbold_0,\,h^t \subseteq \hbold_0$ then $F^r(g^t,h^t) = 1$.
Therefore for any witness $(\gbold,\hbold)$ of $r$\hspace{1mm}, we have 
$g^t \subseteq \gbold,\,h^t \subseteq \hbold$.
\item There are 3 cases :
\begin{enumerate}
\item $g^t \neq g^{t'}$ or $h^t \neq h^{t'}$. Then, since all those functions have the same domain,
we get that $r,t'$ cannot have a common witness  $(\gbold,\hbold)$
because by (1) we must have $g^t \subseteq \gbold,\,h^t \subseteq
\hbold$.
\item $F^t \neq F^{t'}$. Then, since \mbox{$\Dom(F^t)= \Lambda^* \subseteq \Lambda^{s^+}=
\Dom(F^r)$} and $(r,t) \in T$ we know that $F^t \subseteq F^r$.
Since $F^t \neq F^{t'}$ and \\
$\Dom(F^t) =\Dom(F^{t'})$, we get that $F^r$ and $F^{t'}$ aren't
compatible (and therefore there is no common witness)
\item $z^t \neq z^{t'}$. By the previous cases we may assume
that\\
\mbox{$F^t = F^{t'},\,g^t = g^{t'},\,h^t = h^{t'}$} therefore
$z^t,z^{t'}$ are on the same level (See \ref{Sec3-Def-Parameter}
2(d)). We can also see that $z^r$ must be on a greater level
(Remember that $F^t \subseteq F^r$ and $F^r(g^t,h^t) = 1$). Since
$(r,t) \in T$, $z^t,z^r$ are comparable in $\leq^\niceT
\Rightarrow z^{t'},z^r$ are not $\Rightarrow (r,t') \notin T$
\end{enumerate}
\hfill$\square_{\ref{Sec3-Claim-unique-projection}}$
\end{enumerate}

\begin{claim}
\label{Sec3-Claim-exists-projection} For every $i \in A$ there is
$r \in W_{s^+}$ such that $(r,t_i) \in T$
\end{claim}
\par \noindent
Proof:\\
Since $(c_s,c_{s^+}) \in \GG_{s,s^+}$ and this group is generated
by \mbox{$\{(x_{t},x_{t'}):\, (t,t') \in T \cap (J_s \times
J_{s^+})\}$}, there are representations(not necessarily
reduced)\\
 $c_{s_i} = x_{w_1}+\dots+x_{w_n}\,\,c_{s^+} =
x_{r_1}+\dots+x_{r_n}$ such that $(r_n,w_n) \in T$.\\
We may assume that if $1 \leq \ell_1 < \ell_2 \leq n $, then either $r_{\ell_1} \neq r_{\ell_2}$
or $w_{\ell_1} \neq w_{\ell_2}$. 
(Otherwise, we can reduce both representations - remember 
that in those groups $2x =0$). Since $x_{t_i}$ appears in the reduced representation of
$c_{s_i}$, $t_i$ must appear among the $w$-ies. Let $\ell$ be such that $w_\ell = t_i$. 
Now we show that if  $\ell_1 \neq \ell$, then $r_{\ell_1} \neq r_\ell$. 
Assume toward contradiction that $r_{\ell_1} = r_\ell$. By our assumption, 
$w_{\ell_1} \neq w_\ell$. 
Now, we have:
\begin{enumerate}
\item $(r_{\ell_1},w_{\ell_1}),(r_\ell,w_\ell) \in T$
\item $w_\ell \in W_{s_i}$
\item $w_\ell \neq w_{\ell_1}$
\end{enumerate}
this contradicts \ref{Sec3-Claim-unique-projection}.\\
We got that for every $\ell_1 \neq \ell, \,\,\, r_{\ell_1} \neq r_\ell$. This implies that $x_{r_\ell}$
does not cancel, so $r_\ell \in W_{s^+}$ and we are
done.\hfill$\square_{\ref{Sec3-Claim-exists-projection}}$\\

Now choose for each $i \in A \quad r_i \in W_{s^+}$ such that
$(r_i,t_i) \in T$.
\begin{claim}
\label{Sec3-claim-r-ies-distinct} $i < j \Rightarrow r_i \neq r_j$
\end{claim}
\par \noindent
Proof:\\
If we are in case (*1): 
$\{z^{t_i}: \, i \in A\}$  is an antichain. So, $z^{t_i},z^{t_j}$ are not
comparable. Since $z^{r_i} \geq^\niceT z^{t_i}$ and $z^{r_j} \geq^\niceT z^{t_j}$ (See the proof of
\ref{Sec3-Claim-unique-projection} - $z^{r_i},z^{t_i}$ are comparable and $z^{r_i}$ is on greater level), 
We must have $r_i \neq r_j$.\\
If we are in case (*2): assume toward contradiction that $r = r_i
= r_j$. Let $(\gbold,\hbold)$ be a witness for $r$. By
$\ref{Sec3-Claim-unique-projection}\quad g^{t_i},g^{t_j} \subseteq
\gbold,\,\, h^{t_i},h^{t_j} \subseteq \hbold$. Since we are in
case (*2) we get that $\gbold(\alpha_i) = \gbold(\alpha_j)$ but by
the construction \mbox{$\hbold(\alpha_i) < \alpha_j <
\hbold(\alpha_j)$} which contradicts the definition of a witness
(see \ref{Sec3-Def-Parameter} 2(e)).
\hfill$\square_{\ref{Sec3-claim-r-ies-distinct}}$\\
\par \noindent
We got that $W_{s^+}$ is infinite - contradiction. Therefore $M$
must be
rigid.\hfill$\square_{\ref{Sec3-Claim-Not-Iso}}\square_{\ref{Sec3-Theorem-Equiv}}$
\newpage

%
%

\section{ $\lambda > cf(\lambda) > \beth_\omega$}
\noindent
Clearly, for \mbox{$\lambda$ singular  $> \beth_\omega$} we cannot prove
the same result as for\\
 \mbox{$\lambda$ regular $ > \beth_\omega$} (Since
in such game AIS will be able to list all the elements of the two
models). Therefore, we define another type of game.
\begin{definition}
Let $M_1,M_2$ be models with common vocabulary. Let $\niceT$ be a
tree.  We define the game $\Game^*_\niceT(M_1,M_2)$ in the same
way as the definition of $\Game_{\niceT,\mu}$ (See
\ref{Sec1-Def-Game-With-Tree} ) except that in stage $\alpha$ we
demand that the sets $A_1,A_2$ chosen by AIS will satisfy $|A_1
\cup A_2| < 1+ \alpha$ instead of $|A_1 \cup A_2 | < 1 + \mu$. We
say that $M_1,M_2$ are $EF^*_\niceT$ equivalent if ISO has a
winning strategy for $EF^*_\niceT(M_1,M_2)$.
\end{definition}
\begin{remark}
Note that in theorem \ref{Sec2-Theorem-Eq-T,1}, if we replace
$EF_{\niceT,1}$ with $EF^*_\niceT$ we don't get a stronger result,
because for every tree $\niceT$ which satisfies the conditions
there, we can construct another tree $\niceT'$ which satisfies the
conditions, such that $EF_{\niceT',1}$ equivalence would imply
$EF^*_\niceT$ equivalence.
\end{remark}

\begin{theorem}
\label{Sec4-Theorem-EF*-T-Equiv} Suppose that :
\begin{enumerate}
\item $\lambda > \cf(\lambda) = \kappa > \beth_\omega$
\item $\niceT$ is a tree of size $\lambda$ without a $\lambda$
branch
\end{enumerate}
then:\\
There are non-isomorphic models $M_1,M_2$ of size $\lambda$ which
are $ EF^*_{\niceT}$ equivalent.
\end{theorem}
\par \noindent
Proof:\\
Let $\chi$ be a large enough cardinal(for example $\chi =
\beth_7(\lambda)$).
\begin{claim}
\label{Sec4-Claim-exists-inner-model}
We can find $\Mfrak$ such that:
\begin{enumerate}
\item $\Mfrak$ is elementary sub-model of $\niceH(\chi)$
\item $ \lambda+1 \subseteq \Mfrak$
\item for every $ \langle\,(x_i,z_i):\,i < \kappa \rangle$ such that $x_i \in \Mfrak$, $z_i \in \niceT$ for every 
$i < \lambda$  there exists an increasing sequence $\langle i_n :\, n < \omega \rangle $ such that :
      \begin {enumerate}
      \item $\langle (x_{i_n},z_{i_n}):\,n < \omega \rangle \in \Mfrak$ 
      \item if in addition, for every $ \alpha < \lambda $ there is $ i < \kappa$ such that the level of $z_i$ is
      greater then $\alpha$, then we can also have that  $\langle z_{i_n}:\,n < \omega \rangle $ is 
      an antichain in   $\leq^\niceT$
      \end{enumerate}
\end{enumerate}
\end{claim}
\par \noindent
Proof:\\ 
The same proof as the proof of  \ref{Sec3-Claim-exists-inner-model} 
( We are using the fact that $\kappa$ is regular and $\kappa > \beth_\omega$ )
\hfill$\square_{\ref{Sec4-Claim-exists-inner-model}}$ \\

\par \noindent 
Let $\Mfrak$ be as in claim \ref{Sec4-Claim-exists-inner-model}.
Let $\langle \mu_i \, : i < \kappa \rangle$ be an increasing and
continuous sequence such that $\mu_0 = 0$, \mbox{$\mu_i^+ + \aleph_0
< \mu_{i+1} = \cf(\mu_{i+1})$}, $\cup_{i<\kappa} \, \mu_i =
\lambda$ .\\
For every $\alpha < \lambda$ there is a unique $i <
\kappa$, such that $\alpha \in [\mu_i,\mu_{i+1})$. We denote this $i$ by
${\mathbf i}(\alpha)$. \\
\par\noindent
We define a structure parameter $\xfrak$ in the following way:
\begin{definition}
\label{Sec4-Def-Parameter}
\begin{enumerate}

\item $I$ consists of the objects of the form $(u,\Lambda)$ where:
\begin{enumerate}
\item $u \in \lambda^{<\aleph_0}$
\item $\Lambda \in \Mfrak$, $| \Lambda | \leq \aleph_0$, $\Lambda$ is a set of
partial functions with finite domain, from $\lambda$ to $\lambda$.
\end{enumerate}
for $s = (u,\Lambda)$ we denote $u = u^s,\, \Lambda = \Lambda^s$\\
We define $\Gamma(s) = u^s \cup \,\bigcup\{ \Dom(f):\, f \in
\Lambda^s\}$. Note that this a countable set.

\item For $s = (u,\Lambda) \in I$, $J_s$ consists of the objects of the form \mbox{$t
= (u,\Lambda,g,h,F,z)$} where:
\begin{enumerate}
\item $g,h$ are functions from $u$ to $\lambda$
\item $F$ is a function from $\Lambda^2$ to $\{0,1\}$
\item $z \in \niceT$
\item Let $\alpha$ be the level of $z$ in the tree $\niceT$. then
$\alpha$ is minimal under the condition that $\alpha \geq
\mu_{\ibold(x)}$ for every $x$ such that: \\
$h(x) = \mu_{\ibold(x)+1}$ or there are $f_1,f_2 \in \Lambda$ such
that $F(f_1,f_2) =1$ and $f_2(x) = \mu_{\ibold(x)+1}$
\item There is a witness $(\gbold,\hbold)$ for $t$, which means
that :
\begin{enumerate}
\item $\Dom(\gbold) = \Dom(\hbold) \subseteq
\lambda,\,\,\Rang(\gbold) \cup \Rang(\hbold) \subseteq \lambda$
\item $\Gamma(s) \subseteq \Dom(\gbold)$

\item $g \subseteq \gbold, h \subseteq \hbold$
\item for every $(f_1,f_2) \in \Lambda^2$\\
$F(f_1,f_2) = 1$ iff \hspace{1mm} $f_1 \subseteq \gbold \wedge\,\,
f_2 \subseteq \hbold$
\item $\gbold,\hbold$ \hspace{1mm} are weakly increasing
\item $\hbold(x) > x$
\item $\gbold(x) = \gbold(y) \Rightarrow \hbold(x) = \hbold(y)$
\item $\gbold(x) \in [\mu_{\ibold(x)},\mu_{\ibold(x)}^+]$
\item $\hbold(x) \in [\mu_{\ibold(x)},\mu_{\ibold(x)+1}]$
\end{enumerate}
\end{enumerate}

\item $S = I^2$

\item $T$ consists of the pairs $(t_1,t_2) \in J^2$ where :
\begin{enumerate}
\item $t_1,t_2$ have a common witness
\item $z^{t_1},z^{t_2}$ are comparable in the order $\leq^\niceT$
\end{enumerate}

\end{enumerate}
\end{definition}
\begin{fact}
if:
\begin{enumerate}
\item $ s \in I, z \in\niceT$
\item $\gbold,\hbold$ satisfy (i)-(v)
\item $\bigcup\{\mu_{\ibold(x)} :\, \hbold(x) = \mu_{\ibold(x)+1}\}  \subset \alpha$ where $\alpha$ is the level of $z$
\end{enumerate}
then:
\begin{enumerate}
\item there is unique $t \in J_s$ such that $(\gbold,\hbold)$ is a
witness for t, and $z^t \leq^\niceT z$. we denote $t =
t(s,\gbold,\hbold,z)$
\item if :
\begin{enumerate}
\item $\gbold',\hbold',z'$ satisfy the conditions in (1)
\item $z,z'$ are comparable in $\leq^\niceT$
\item  $\gbold',\hbold'$ are compatible with $\gbold,\hbold$ respectively
\end{enumerate}
then :\\
$t(s,\gbold,\hbold,z) = t(s,\gbold',\hbold',z')$ \\
\end{enumerate}
\end{fact}
\noindent
Let $M = M_{\mathfrak x}$ be the corresponding model. We can check
that $||M|| = \lambda$ . Let $a_* =
0_{\GG_{(\emptyset,\emptyset)}},\quad b_* =
x_{(\emptyset,\emptyset,\emptyset,\emptyset,\emptyset,z_*)}$ where $z_*$ is the
root of $\niceT$ (without loss of generality there is a root).\\
Define $M_1 = (M,a_*), M_2 = (M,b_*)$.
\begin{claim}
\label{Sec4-Claim-Equiv} $M_1,M_2$ are $EF^*_{\niceT}$ equivalent.
\end{claim}
\par \noindent
We describe a winning strategy for ISO - this is very similar to
the proof of  \ref{Sec2-Claim-Equiv-Models}, so we omit the
details. We use the definitions in \ref{Sec2-Def-W} . In every stage $\alpha$ of the game, ISO will choose
a function $\gbold_\alpha$, such that :
\begin{enumerate}
\item $\gbold_0 = \emptyset$
\item $\gbold_\alpha \in \niceW_{\ibold(\alpha)+1}$
\item $\beta < \alpha \Rightarrow \gbold_\beta \leq^\niceW \gbold_\alpha$
\item if in stage $\alpha$ AIS chose the sets $A_1,A_2$ then\\
for each $s \in I$, if $\GG_s \cap (A_1 \cup A_2) \neq \emptyset$
then $\Gamma(s) \subseteq \Dom(\gbold_\alpha)$
\end{enumerate}
Now if $\alpha = \beta + 1$ and in stage $\alpha$ AIS chose the
sets $A_1,A_2$ and the node $z_\alpha$, ISO will define
$\hbold_\alpha = h_{\gbold_\alpha}$, and then define $f_\alpha$ by
:
\begin{enumerate}
\item $\Dom(f_\alpha) = \bigcup\{\GG_s : \Gamma(s) \subseteq
\Dom(\gbold_\alpha)\}$
\item for each $s$ such that $\GG_s \subseteq \Dom(f_\alpha)$,\\
$f_\alpha(0_{\GG_s}) = x_t$ where $t =
t(s,\gbold_\alpha,\hbold_\alpha,z_\alpha)$\hfill$\square_{\ref{Sec4-Claim-Equiv}}$
\end{enumerate}

\begin{claim}
\label{Sec4-Claim-Not-Isomorphic}
$M_1,M_2$ are not isomorphic
\end{claim}
\par \noindent
Proof:\\
It is enough to show that $M$ is rigid. The proof is very similar
to the proof of \ref{Sec3-Claim-Not-Iso} . Assume toward
contradiction that $f \neq id$ is an automorphism of $M$. Denote
$W_s = \{t \in J_s :\, x_t\,\, \textrm{is in the reduced
representation of}\,\, c_s \}$. Since $f \neq id$ there is  $s^* =
(u^*,\Lambda^*)$ such that $W_{s^*} \neq \emptyset$.
\par \noindent
Case (*1):\\
We can find $\langle s_\theta,t_\theta,\alpha_\theta :
\theta < \kappa \rangle$ such that:
\begin{enumerate}
\item $ s_\theta \in J, s_\theta = (u^* \cup\ \{\alpha_\theta\}, \Lambda^*)$
\item $ t_\theta \in W_{s_\theta}$
\item $ h^{t_\theta}(\alpha_\theta) = \mu_{\ibold(\alpha_\theta)+1}$
\item $\theta < \varepsilon < \kappa \Rightarrow \ibold(\alpha_\theta) <
\ibold(\alpha_{\varepsilon})$
\end{enumerate}
In this case, note that the level of $z^{t_\theta}$ must be $\geq
\mu_{\ibold(\alpha_\theta)}$.
\par \noindent
Case (*2):\\
We cannot find such a sequence. Therefore, for every
large enough $i < \kappa$, for every $\alpha$ such that
$\ibold(\alpha) = i $, for $s(\alpha) = (u^* \cup
\{\alpha\},\Lambda^*)$, for every $t \in W_{s(\alpha)}$,
\hspace{1mm}$h^t(\alpha) < \mu_{i +1}$.\\
Choose $i^*$ which satisfies this and $\mu_{i^*} > \mu$.\\
We can find $\langle t_\theta,s_\theta,\alpha_\theta : \theta <
\mu_{i^*+1} \rangle$ such that :
\begin{enumerate}
\item $s_\theta \in I,\,t_\theta \in W_{s_\theta}$
\item $\ibold(\alpha_\theta) = i^*$
\item $\theta < \varepsilon \Rightarrow h^{t_\theta}(\alpha_\theta)
< \alpha_{\varepsilon}\,\, (<
h^{t_\varepsilon}(\alpha_\varepsilon)\,)$
\end{enumerate}
Since $\mu_{i^*+1} = cf(\mu_{i^*+1}) >  \mu_{i^*}^+$ and for every
$\theta$ we have $g^t_\theta(x) \leq \mu_{i^*}^+$ (This is by
\ref{Sec4-Def-Parameter}(2)(e)(viii)), we may assume that
$g^{t_\theta}(\alpha_\theta)$ is constant.\\
\par \noindent
Now, in both cases, we proceed in a similar way to the proof of \ref{Sec3-Claim-Not-Iso}. Using \ref{Sec4-Claim-exists-inner-model}, we choose $A \subset \kappa $ such that:
\begin{enumerate} 
\item $|A| = \aleph_0$
\item $\langle W_{s_\theta}:\, \theta \in A \rangle \in \Mfrak$
\item if we are in case (*1) then  $\langle z^{t_\theta}:\, \theta \in A \rangle $ is an antichain in $\leq^\niceT$
      ( We can demand this because in case (*1) the levels of the $z^{t_\theta}$-ies aren't 
      bounded in $\lambda$ - See    \ref{Sec4-Claim-exists-inner-model} )
\end{enumerate}
Define $s^+ \in I$ by \mbox{$s^+ = (\emptyset,\Lambda^* \cup \{g^t,h^t:\,t \in
W_{s_\theta},\, \theta \in A\})$}.
\begin{claim}
\label{Sec4-Claim-unique-projection} For every $\theta \in
A$,\hspace{1mm} if : $r \in J_{s^+},\,t \in W_{s_\theta}, (r,t)
\in
T$\\
then :
\begin{enumerate}
\item if $(\gbold,\hbold)$ is a witness for $r$ then $g^t
\subseteq \gbold,\,h^t \subseteq \hbold$
\item if $t \neq t' \in J_{s_\theta} $ then $(r,t') \notin T$
\end{enumerate}
\end{claim}
\par \noindent
Proof: see the proof of \ref{Sec3-Claim-unique-projection}. 
\hfill$\square_{\ref{Sec4-Claim-unique-projection}}$
\begin{claim}
\label{Sec4-Claim-exists-projection} For every $\theta \in A$
there is $r \in W_{s^+}$ such that $(r,t_\theta) \in T$
\end{claim}
\par \noindent
Proof: see the proof of \ref{Sec3-Claim-exists-projection}
\hfill$\square_{\ref{Sec4-Claim-exists-projection}}$\\
\par \noindent
Now, using \ref{Sec4-Claim-exists-projection}, we choose for each $\theta \in A$, $r_\theta \in W_{s^+}$, such
that $(t_{\theta},r_{\theta}) \in T$.
\begin{claim}
\label{Sec4-claim-r-ies-distinct} $\theta < \varepsilon
\Rightarrow r_\theta \neq r_\varepsilon$
\end{claim}
\par \noindent
Proof:\\
If we are in case (*1):\\
$z^{t_\theta},z^{t_\varepsilon}$
are not comparable. But, $z^{r_\theta} \geq^\niceT z^{t_\theta}$ because they are comparable and 
$z^{r_\theta}$ is on greater level, since that level is determined by \ref{Sec4-Def-Parameter} 2(d).
By the same argument, $z^{r_\varepsilon} \geq^\niceT z^{t_\varepsilon}$. Therefore, $z^{r_\varepsilon},\,z^{r_\theta}$
aren't comparable, so $r_\theta \neq r_\varepsilon$.\\
If we are in case (*2):\\
Assume toward contradiction that $r =r_\theta = r_\varepsilon$.\\
Let $(\gbold,\hbold)$ be a witness for $r$. By $\ref{Sec4-Claim-unique-projection}\quad
g^{t_\theta},g^{t_\varepsilon} \subseteq \gbold,\,\,
h^{t_\theta},h^{t_\varepsilon} \subseteq \hbold$. Since we are in
case (*2) we get that $\gbold(\alpha_\theta) =
\gbold(\alpha_\varepsilon)$ and \mbox{$\hbold(\alpha_\theta) <
\alpha_\varepsilon < \hbold(\alpha_\varepsilon)$} \\
which contradicts the definition of a witness (see
\ref{Sec4-Def-Parameter} 2(e)).
\hfill$\square_{\ref{Sec4-claim-r-ies-distinct}}$\\
\par \noindent
We got that $W_{s^+}$ is infinite - contradiction. Therefore, $M$
must be rigid. \hfill$\square_{\ref{Sec4-Claim-Not-Isomorphic}}\square_{\ref{Sec4-Theorem-EF*-T-Equiv}}$
\newpage
\bibliographystyle{alpha}

\begin{thebibliography}{0}

\bibitem {Chang} C.C Chang, Some remarks on the model theory of infinitary languages, in \emph{The syntax and semantics of infinitary languages}, Lecture notes in Mathematics, 72, J. Barwise ed. (Springer, Berlin, 1968) pp. 36-63

\bibitem{Hyttinen} T.Hyttinen and H. Tuuri, Constructing strongly equivalent nonisomorphic models for unstable theories, Annals of pure And Applied Logic 52(1991) pp.203-248

\bibitem{Hodges}  Hodges, Wilfrid Model theory. Encyclopedia of Mathematics and its Applications, 42. Cambridge University Press,Cambridge, 1993.  ISBN: 0-521-30442-3 (Reviewer: J. M. Plotkin)

\bibitem{SH-836} S.Shelah A long EF equivalence non isomorphic models - to appear (SH 836 in Shelah archive)

\bibitem{SH-RGCH} S.Shelah, The generalized  continum hypothesis revisited, Israel J. Math 116 (2000) pp. 285-321

\bibitem{SH-Lambda-Sing-Non-Iso} S.Shelah, Existence of many $L_{\infty,\lambda}$ non isomorphic models of power $\lambda$ for $\lambda$ singular with $\lambda^\omega = \lambda$, Notre Dame J. Formal Logic 25(1984) pp. 97-104

\bibitem{SH-Many-Models-For-Unstable} S.Shelah, Existence of many $L_{\infty,\lambda}$ non isomorphic models of $T$ of power $\lambda$, Annals of Pure and Applied Logic 34 (1987) pp. 291-310

\bibitem{SH-Classification} S.Shelah, Classification Theory, Stud. Logic Found. Math. 92 (North Holland, Amsterdam, 2nd rev. ed. 1990)

\bibitem{Scott} D.S. Scott Logic with denumerably long formulas and finite strings of quantifiers in 
{\emph{The theory of models}}, J.W.Adisson, L.Henkin and A.Tarsky , eds. (North holland, Amsterdam, 1965) pp. 329-341 


\end{thebibliography}

\end{document}